\documentclass[12pt]{article}
\usepackage{amsmath, amssymb}

\newcommand{\qed}{$\hfill \Box$}
\newcommand{\proof}{\medskip \noindent {\bf Proof:}\qquad}

\newtheorem{theorem}{Theorem}[section]

\newtheorem{remark}{Remark}[section]

\newlength{\IndentI}
\newlength{\IndentII}
\newlength{\IndentIII}
\newlength{\WidthI}
\newlength{\WidthII}
\newlength{\WidthIII}
\title{On a simple strategy weakly forcing the strong law of large
  numbers in the bounded forecasting game} \author{
  Masayuki Kumon\\
  Risk Analysis Research Center\\
  Institute of Statistical Mathematics\\
  and\\
  Akimichi Takemura\\
  Graduate School of Information Science and Technology\\
  University of Tokyo }

\date{November, 2005}

\begin{document}
\maketitle

\begin{abstract}
  In the framework of the game-theoretic probability of Shafer and
  Vovk (2001) it is of basic importance to construct an explicit
  strategy weakly forcing the strong law of large numbers (SLLN) in
  the bounded forecasting game.  We present a simple finite-memory
  strategy based on the past average of Reality's moves, which weakly
  forces the strong law of large numbers with the convergence rate of
  $O(\sqrt{\log n/n})$.  Our proof is very simple compared to a
  corresponding measure-theoretic result of Azuma (1967) on bounded
  martingale differences and this illustrates  effectiveness of
  game-theoretic approach.  We also discuss one-sided protocols and
  extension of results to linear protocols in general dimension.
\end{abstract}


\section{Introduction}
The book by Shafer and Vovk (2001) established the whole new field of
game-theoretic probability and finance.  Their framework provides an
attractive alternative foundation of probability theory.  Compared to
the conventional measure theoretic probability, the game theoretic
probability treats the sets of measure zero in a very explicit way
when proving various probabilistic laws, such as the strong law of
large numbers. In a game-theoretic proof, we can explicitly
describe the behavior of the paths on a set of measure zero, whereas
in measure-theoretic proofs the sets of measure zero are often simply
ignored. This feature of game-theoretic probability is well
illustrated in the explicit construction of Skeptic's strategy forcing
SLLN in Chapter 3 of Shafer and Vovk (2001).

However the strategy given in Chapter 3 of Shafer and Vovk (2001),
which we call a mixture $\epsilon$-strategy in this paper, is not yet
satisfactory, in the sense that it requires combination of infinite
number of ``accounts'' and it needs to keep all the past moves of
Reality in memory.  We summarize their construction in Appendix in a
somewhat more general form than in Chapter 3 of Shafer and Vovk
(2001).  In fact in Section 3.5 of their book, Shafer and Vovk pose
the question of required memory for strategies forcing SLLN.

In this paper we prove that a very simple single strategy, based only
on the past average of Reality's moves is weakly forcing SLLN.
Furthermore it weakly forces SLLN with the convergence rate of
$O(\sqrt{\log n/n})$.  In this sense, our result is a substantial
improvement over the mixture $\epsilon$-strategy of Shafer and Vovk.
Since $\epsilon$-strategies are used as essential building blocks for
the ``defensive forecasting'' (\cite{vovk-takemura-shafer-defensive}),
the performance of defensive forecasting might be improved by
incorporating our simple strategy.

Our thinking was very much influenced by the detailed analysis by Takeuchi
(\cite{takeuchi:2004b} and Chapter 5 of \cite{takeuchi:2004}) of the
optimum strategy of Skeptic in the games, which are favorable for
Skeptic.  We should also mention that the intuition behind our
strategy is already discussed several times throughout the book by
Shafer and Vovk (see e.g.\ Section 5.2).  Our contribution is in
proving that the strategy based on the past average of Reality's moves
is actually weakly forcing SLLN.

In this paper we only consider weakly forcing by a strategy.  A strategy
weakly forcing an event $E$ can be transformed to a strategy forcing
$E$ as in Lemma 3.1 of Shafer and Vovk (2001).  We do not present
anything new for this step of the argument.

The organization of this paper is as follows.  In Section
\ref{sec:motivation} we formulate the bounded forecasting game and
motivate the strategy based on the past average of Reality's moves as
an approximately optimum $\epsilon$-strategy.  In Section
\ref{sec:slln} we prove that our strategy is weakly forcing SLLN with
the convergence rate of $O(\sqrt{\log n/n})$.  In Section
\ref{sec:one-sided} we consider one-sided protocol and prove that
one-sided version of our strategy weakly forces one-sided SLLN with
the same order. In Section \ref{sec:linear-protocol} we treat a
multivariate extension to linear protocols.  In Appendix we give a
summary of the mixture $\epsilon$-strategy in Chapter 3 of Shafer and
Vovk (2001).

\section{Approximately optimum single $\epsilon$-strategy for the
  bounded forecasting game}
\label{sec:motivation}

Consider the bounded forecasting game in Section 3.2 of Shafer and
Vovk (2001).

\medskip
\noindent
\textsc{Bounded Forecasting Game}\\
\textbf{Protocol:}

\parshape=6
\IndentI   \WidthI
\IndentI   \WidthI
\IndentII  \WidthII
\IndentII  \WidthII
\IndentII  \WidthII
\IndentI   \WidthI
\noindent
${\cal K}_0 :=1$.\\
FOR  $n=1, 2, \dots$:\\
  Skeptic announces $M_n\in{\mathbb R}$.\\
  Reality announces $x_n\in [-1,1]$.\\
  ${\cal K}_n := {\cal K}_{n-1} + M_n x_n$.\\
END FOR

\medskip
\noindent
For a fixed $\epsilon$, $|\epsilon| < 1$, the $\epsilon$-strategy sets
$ M_n = \epsilon  {\cal K}_{n-1}$.   Under this strategy Skeptic's capital
process ${\cal K}_n$ is written as 
${\cal K}_n = \prod_{i=1}^n (1+\epsilon x_i)$ or
\[
\log {\cal K}_n = \sum_{i=1}^n \log (1+\epsilon x_i).
\]
For sufficiently small $|\epsilon|$, $\log {\cal K}_n$ is approximated as
\[
\log {\cal K}_n \simeq \epsilon \sum_{i=1}^n x_i - \frac{1}{2} \epsilon^2
\sum_{i=1}^n x_i^2.
\]
The right-hand side is maximized by taking
\[
\epsilon = \frac{\sum_{i=1}^n x_i}{\sum_{i=1}^n x_i^2}.
\]

In particular in the fair-coin game, where $x_n$ is restricted as $x_n =
\pm 1$, approximately optimum $\epsilon$ is given as
\[
\epsilon = \bar x_n = \frac{1}{n} \sum_{i=1}^n x_i .
\]
Actually, as shown by Takeuchi (\cite{takeuchi:2004b}), it is easy to check
that $\epsilon = \bar x_n$ exactly maximizes $\prod_{i=1}^n (1+\epsilon
x_i)$ for the case of the fair-coin game.  Recently Kumon, Takemura
and Takeuchi (2005) \cite{ktt1} give a detailed analysis of Bayesian
strategies for the biased-coin games, which include the strategy
$\epsilon = \bar x_n$ as a special case.

Of course, the above approximately optimum $\epsilon$ is chosen in
hindsight, i.e., we can choose optimum $\epsilon$ after seeing the moves
$x_1,\dots, x_n$.  However it suggests to choose $M_n$ based on the past
average $\bar x_{n-1}$ of Reality's moves.  Therefore consider a
strategy ${\cal P}={\cal P}^c$ 
\begin{equation}
\label{eq:our-strategy}
M_n = c \bar x_{n-1} {\cal K}_{n-1}.
\end{equation}
In the next section we prove that for $0< c \le 1/2$ this strategy is
weakly forcing SLLN.  The restriction $0< c \le 1/2$ is just for
convenience for the proof and 
in Kumon, Takemura and Takeuchi (2005) we consider 
$c=1$ for biased-coin games.  


Compared to a single fixed $\epsilon$-strategy $M_n = \epsilon {\cal K}_{n-1}$
or the mixture $\epsilon$-strategy in Chapter 3 of Shafer and
Vovk (2001), letting $\epsilon = c\bar x_{n-1}$ depend on $\bar x_{n-1}$
seems to be reasonable from the viewpoint of effectiveness of Skeptic's
strategy.  The basic reason is that as $\bar x_{n-1}$ deviates more from
the origin, Skeptic should try to exploit this bias in Reality's moves by
betting a larger amount.  Clearly this reasoning is shaky because for each
round Skeptic has to move first and Reality can decide her move after
seeing Skeptic's move.  However in the next section we show that the
strategy  in (\ref{eq:our-strategy}) is indeed weakly forcing SLLN with
the convergence rate of $O(\sqrt{\log n/n})$.

\section{Weakly forcing SLLN  by past averages}
\label{sec:slln}


In this section we prove the following result.

\begin{theorem}  \label{thm:main} 
  In the bounded forecasting game, if Skeptic uses the
  strategy (\ref{eq:our-strategy}) with \qquad $0< c \le 1/2$, then
  $\limsup_n {\cal K}_n =\infty$ for each path $\xi=x_1 x_2 \dots$ of
  Reality's moves such that
\begin{equation}
\label{eq:thm}
\limsup_n \frac{\sqrt{n}|\bar x_n|}{\sqrt{\log n}} > 1.
\end{equation}
\end{theorem}

This theorem states that the strategy (\ref{eq:our-strategy}) weakly
forces that $\bar x_n$ converges to 0 with the convergence rate of 
$O(\sqrt{\log n/n})$.  Therefore it is much stronger than the mixture
$\epsilon$-strategy in Chapter 3 of Shafer and Vovk (2001), which only
forces convergence to 0.  A corresponding measure theoretic result was
stated in Theorem 1 of  Azuma (1967) as discussed in Remark
\ref{rem:3.1} at the end of this section.
The rest of this section is devoted to a proof of Theorem \ref{thm:main}.

By comparing $1,1/2,1/3,\dots,$ and the integral of  $1/x$  we have
\[
\log(n+1)=\int_1^{n+1} \frac{1}{x} dx
\le               1 + \frac{1}{2} + \frac{1}{3} + \dots + \frac{1}{n}
\le 1 + \int_1^n \frac{1}{x}dx = 1 + \log n.
\]

Next, by summing up the term
\[
\frac{1}{(k-1)k} = \frac{1}{k-1} - \frac{1}{k}
\]
from $k=i$ to $n$, we have
\begin{equation}
\label{eq:wabun0}
\sum_{k=i}^n \frac{1}{(k-1)k} = \frac{1}{i-1}- \frac{1}{n}
\end{equation}
or 
\[
\frac{1}{i-1} = \sum_{k=i}^n \frac{1}{(k-1)k} + \frac{1}{n}.
\]

Now consider the sum
\begin{equation}
\label{eq:sum-x-bar}
\sum_{i=2}^n \frac{i}{i-1} \bar  x_i^2 
=\sum_{i=2}^n \frac{1}{(i-1)i} (x_1 + \dots+x_i)^2 .
\end{equation}
When we expand the right-hand side, the coefficient of the term $x_j
x_k$, $j<k$, is given by
\[
2\times \left(
\frac{1}{(k-1)k} +\frac{1}{k(k+1)} + \dots + \frac{1}{(n-1)n} \right)
= 2\times \left(\frac{1}{k-1} - \frac{1}{n} \right).
\]
We now consider the coefficient of  $x_j^2$.
In (\ref{eq:sum-x-bar}) we have the sum from $i=2$ and we need to treat 
$x_1$ separately.  The coefficient of $x_1^2$ is $1-1/n$ from
(\ref{eq:wabun0}).  For $j\ge 2$, the coefficient of $x_j^2$ 
is given by $1/(j-1)-1/n$, as in the case of the cross terms.  Therefore
\begin{eqnarray}
\label{eq:sum-x-bar1}
\sum_{i=2}^n \frac{i}{i-1} \bar  x_i^2 
&=& 2 \sum_{1\le j < i \le n}\Bigl(\frac{1}{i-1}-\frac{1}{n}\Bigr)  
x_j x_i +\Bigl(1-\frac{1}{n}\Bigr) x_1^2 + 
\sum_{i=2}^n \Bigl(\frac{1}{i-1}-\frac{1}{n}\Bigr) x_i^2
\nonumber \\
&=& 2 \sum_{1\le j < i \le n} \frac{1}{i-1} x_j x_i 
- \frac{2}{n}\sum_{1\le j < i \le n} x_j x_i + x_1^2+ 
\sum_{i=2}^n \frac{1}{i-1} x_i^2
-\frac{1}{n} \sum_{i=1}^n x_i^2
\nonumber \\
&=& 2 \sum_{i=2}^n \bar x_{i-1} x_i - n \bar x_n^2 
+ x_1^2 + \sum_{i=2}^n \frac{1}{i-1} x_i^2.
\end{eqnarray}
Write $\bar x_0=0$. Then the first sum on the right-hand side can be
written from $i=1$. Under this notational convention
(\ref{eq:sum-x-bar1}) is rewritten as
\begin{equation}
\label{eq:sum-x-bar2}
\sum_{i=1}^n \bar x_{i-1} x_i 
= \frac{1}{2} \sum_{i=2}^n \frac{i}{i-1} \bar  x_i^2  + \frac{n}{2} \bar x_n^2
- \frac{1}{2}\Bigl(x_1^2 + \sum_{i=2}^n \frac{1}{i-1} x_i^2\Bigr).
\end{equation}
\bigskip

Now the capital process ${\cal K}_n={\cal K}_n^{\cal P}$ of (\ref{eq:our-strategy}) 
is written as
\[
{\cal K}_n = \prod_{i=1}^n (1 + c \bar x_{i-1} x_i) .
\]
As in Chapter 3 of Shafer and Vovk (2001) we use 
\[
\log(1+t) \ge t - t^2, \qquad |t|\le 1/2.
\]
Then for $0< c \le 1/2$  we have
\begin{eqnarray*}
\log {\cal K}_n
&=&
 \sum_{i=1}^n \log (1 + c \bar x_{i-1} x_i)
\nonumber \\
&\ge& c \sum_{i=1}^n \bar x_{i-1} x_i - c^2 \sum_{i=1}^n \bar
x_{i-1}^2 x_i^2,
\end{eqnarray*}
and under the restriction $|x_n|\le 1, \forall n$, we can further bound
$\log {\cal K}_n$ from below as
\begin{equation}
\label{eq:logkn}
\log {\cal K}_n \ge 
c \sum_{i=1}^n \bar x_{i-1} x_i - c^2 \sum_{i=1}^n \bar
x_{i-1}^2 .
\end{equation}

By considering the restriction $|x_n|\le 1$, (\ref{eq:sum-x-bar2}) 
is bounded from below as
\begin{eqnarray*}
\sum_{i=1}^n \bar x_{i-1} x_i 
&\ge &
\frac{1}{2} \sum_{i=2}^n \frac{i}{i-1} \bar  x_i^2  + \frac{n}{2} \bar x_n^2
- \frac{1}{2}\Bigl(1 + \sum_{i=2}^n \frac{1}{i-1}\Bigr)
\\
&\ge &
\frac{1}{2} \sum_{i=2}^n \frac{i}{i-1} \bar  x_i^2  + \frac{n}{2} \bar x_n^2
- \frac{1}{2}(2 + \log(n-1))
\\
&\ge &
\frac{1}{2} \sum_{i=2}^n \bar  x_i^2  + \frac{n}{2} \bar x_n^2
- \frac{1}{2}(2 + \log(n-1))
\\
&\ge &
\frac{1}{2} \sum_{i=1}^n \bar  x_i^2  + \frac{n}{2} \bar x_n^2
- \frac{1}{2}(3 + \log(n-1)) .
\end{eqnarray*}
Since  $c \le 1/2$, substituting this into
(\ref{eq:logkn}) yields
\begin{eqnarray*}
\log {\cal K}_n &\ge&
c\Bigl(\frac{1}{2} - c\Bigr) \sum_{i=1}^n \bar x_{i-1}^2 + 
 c\frac{n}{2} \bar x_n^2 - \frac{c}{2}(3 + \log(n-1)) 
\\
&\ge& c \frac{n}{2} \bar x_n^2 - \frac{c}{2}(3 + \log(n-1))
\\
&\ge&
\frac{c}{2} (n\bar x_n^2 - \log n) - \frac{3}{2} c
\\
&=& \frac{c}{2} \log n \left( \frac{n \bar x_n^2}{\log n} -1 \right) - \frac{3}{2}c.
\end{eqnarray*}

Now if $\limsup_n \sqrt{n}|\bar x_n|/\sqrt{\log n} > 1$, then $\limsup_n
\log {\cal K}_n = +\infty$  because $\log n \uparrow\infty$.   This proves
the theorem.

\begin{remark}
\label{rem:3.1}
In the framework of the conventional measure theoretic probability,
a strong law of large numbers analogous to Theorem 3.1 can be proved
using Azuma-Hoeffding-Bennett inequality (Appendix A.7 of Vovk,
Gammerman and Shafer (2005), Section 2.4 of Dembo and Zeitouni
(1998), Azuma (1967), Hoeffding (1963), Bennett (1962)).  
Let $X_1, X_2, \dots$ be a sequence of martingale differences such that
$|X_n|\le 1$, $\forall n$.  Then for any $\epsilon >0$
\[
P( |\bar X_n| \ge \epsilon) \le 2 \exp( -n\epsilon^2/2).
\]
Fix an arbitrary $\alpha > 1/2$.  Then for any $\epsilon >0$
\[
\sum_n P(|\bar X_n| \ge \epsilon (\log n)^\alpha /\sqrt{n})
\le \sum_n  \exp (-\frac{\epsilon^2}{2}(\log n)^{2\alpha}) \ < \infty.
\]
Therefore by Borel-Cantelli $\sqrt{n}|\bar X_n|/(\log n)^\alpha
\rightarrow 0$ almost surely.  
Actually Theorem 1 of Azuma (1967) states the following stronger result
\begin{equation}
\label{eq:azuma}
\limsup_{n\rightarrow\infty}
 \frac{\sqrt{n} \bar x_n}{\sqrt{\log n}} \le \sqrt{2} \qquad {\rm a.s.}
\end{equation}
Although our Theorem \ref{thm:main} is better in the constant factor
of $\sqrt{2}$, Azuma's result (\ref{eq:azuma}) and our result
(\ref{eq:thm}) are virtually the same.  However we want to emphasize
that our game theoretic proof requires much less mathematical
background than the measure theoretic proof.  Also see the factor of
$\sqrt{3/2}$ in the one-sided version of our result in Theorem
\ref{thm:one-sided} below.
\end{remark}

\section{One-sided protocol}
\label{sec:one-sided}

In this section we consider one-sided bounded forecasting game where
$M_n$ is restricted to be nonnegative ($M_n \ge 0$), i.e.\ Skeptic is
only allowed to buy tickets.  We also consider the restriction
$M_n \le 0$. In Chapter 3 of Shaver and Vovk, weak forcing of SLLN is
proved by combining positive and negative one-sided strategies,
whereas in the previous section we proved that a single strategy
${\cal P}={\cal P}^c$ weakly forces SLLN.  Therefore it is of interest
to investigate whether one-sided version of our strategy weakly forces
one-sided SLLN.  We adopt the same notations as Section 5 of Kumon,
Takemura and Takeuchi (2005), where one-sided protocols for
biased-coin games are studied.

For the positive one-sided case consider the strategy ${\cal P}^+$ with
\[
 M_n = c\bar x_{n-1}^+ {\cal K}_{n-1}, \qquad 
\bar x_{n-1}^+ = \max(\bar x_{n-1},0).
\]
Similarly we consider negative one-sided strategy ${\cal P}^-$
with $M_n = -c\bar x_{n-1}^- {\cal K}_{n-1}$,  
$\bar x_{n-1}^- = \max(-\bar x_{n-1},0)$.

For these protocols we have the following theorem.

\begin{theorem}  \label{thm:one-sided} 
  If Skeptic uses the
  strategy ${\cal P}^+$ with $0< c \le 1/2$, then
  $\limsup_n {\cal K}_n =\infty$ for each path $\xi=x_1 x_2 \dots$ of
  Reality's moves such that 
\[
\limsup_n \frac{\sqrt{n} \bar x_n}{\sqrt{\log n}} > \sqrt{\frac{3}{2}}.
\]
Similarly if Skeptic uses the
  strategy ${\cal P}^-$ with $0< c \le 1/2$, then
$\limsup_n {\cal K}_n =\infty$ for each path $\xi=x_1 x_2 \dots$ of
  Reality's moves such that 
$\liminf_n  \sqrt{n} \bar x_n/\sqrt{\log n} < -\sqrt{3/2}$.
\end{theorem}


The rest of this section is devoted to a proof of this theorem for
${\cal P}^+$.
If $\bar x_n$ is eventually all nonnegative, then
the behavior of the capital process 
${\cal K}^{\cal P}_n$ and ${\cal K}_n^{{\cal P}^+}$ are asymptotically
equivalent except for a constant factor reflecting some initial segment
of Reality's path $\xi$.  Then the theorem follows from Theorem
\ref{thm:main}.
On the other hand if $\bar x_n$ is eventually negative,
then ${\cal K}_n^{{\cal P}^+}$ stays constant and
Theorem \ref{thm:one-sided}  holds trivially.
Therefore we only need to consider the case that $\bar x_n$ changes sign
infinitely often.
Note that at time $n$ when $\bar x_n $ changes the sign, the overshoot
is bounded as
\[
|\bar x_n| \le  1/n.
\]

We consider capital process after a sufficiently large time $n_0$ such
that $\bar x_{n_0} \simeq 0$, and proceed to divide the sequence
$\{\bar x_n\}$ into the following two types of blocks.  For $n_0 \le k 
\le l - 1$, consider a block $\{k, \dots , l-1\}$.  We call it a
{\it nonnegative block} if
\[ 
\bar x_{k-1} < 0, \bar x_{k} \ge 0, \bar x_{k+1} \ge 0,
\dots ,
\bar x_{l-1} \ge 0, \bar x_{l} < 0. 
\]
Similarly we call it a {\it negative block} if
\[
\bar x_{k-1} \ge 0, \bar x_k < 0, \bar x_{k+1} < 0, \dots
,
\bar x_{l-1} < 0, \bar x_{l} \ge 0.  
\]
By definition, negative and nonnegative blocks are alternating.

For a nonnegative block
\begin{equation}
\label{eq:nonnegative-block}
{ \cal K}^{{\cal P}^+}_{l} = { \cal K}^{{\cal P}^+}_{k} \prod_{i=k+1}^l
(1+ c \bar x_{i-1}x_i)
\end{equation}
whereas for a negative block 
${ \cal K}^{{\cal P}^+}_{l} = { \cal K}^{{\cal P}^+}_{k}$.  Taking the
logarithm of (\ref{eq:nonnegative-block}) we have
\begin{align*}
\log { \cal K}^{{\cal P}^+}_{l}  - \log { \cal K}^{{\cal P}^+}_{k}
&=  \sum_{i=k+1}^l \log (1+c \bar x_{i-1} x_i)\\
&\ge c \sum_{i=k+1}^l \bar x_{i-1} x_i - c^2 \sum_{i=k+1}^l 
    \bar x_{i-1}^2 x_i^2\\
&\ge c \sum_{i=k+1}^l \bar x_{i-1} x_i - c^2 \sum_{i=k+1}^l \bar x_{i-1}^2.
\end{align*}  
From (\ref{eq:sum-x-bar2}) it follows
\begin{align*}
\sum_{i=k+1}^l \bar x_{i-1} x_i 
&= \frac{1}{2} \sum_{i=k+1}^l \frac{i}{i-1} \bar  x_i^2  + \frac{l}{2}
\bar x_l^2 - \frac{k}{2} \bar x_k^2
- \frac{1}{2}\sum_{i=k+1}^l \frac{1}{i-1} x_i^2 \\
&\ge 
\frac{1}{2} \sum_{i=k+1}^l \frac{i}{i-1} \bar  x_i^2  
 - \frac{1}{2k} 
- \frac{1}{2}\left(\frac{1}{k} + \log \frac{l}{k} \right).  
\end{align*}
In the above, we used the approximation formula
\begin{align*}
\frac{1}{m} + \frac{1}{m + 1} + \cdots + \frac{1}{n} \le 
\int_m^n \frac{dx}{x} + \frac{1}{m} = 
\log \frac{n}{m} + \frac{1}{m}.
\end{align*}
Thus we obtain
\begin{equation}
\label{eq:block-end}
\log { \cal K}^{{\cal P}^+}_{l}  - \log { \cal K}^{{\cal P}^+}_{k}
\ge  - \frac{c}{k} 
- \frac{c}{2}\log \frac{l}{k}.
\end{equation}

Now starting at $n_0$, we consider adding the right-hand side of 
(\ref{eq:block-end}) for
nonnegative blocks and 
$0= \log { \cal K}^{{\cal P}^+}_{l}  - \log { \cal K}^{{\cal
    P}^+}_{k}$
for negative blocks.  Then after passing sufficiently many  
blocks, 
$\log { \cal K}^{{\cal P}^+}_{l}  - \log { \cal K}^{{\cal P}^+}_{k}$ 
behave as (\ref{eq:block-end}) 
during half number of the entire blocks. 
Therefore at the beginning $n_k$ of the last nonnegative block, we
have
\begin{equation}
\label{eq:logk-diff}
\log { \cal K}^{{\cal P}^+}_{n_k} - \log { \cal K}^{{\cal P}^+}_{n_0} 
\ge - \frac{c}{2}\log \frac{n_k}{n_0} - \frac{c}{4}\log \frac{n_k}{n_0} 
+o(1) = - \frac{3c}{4}\log \frac{n_k}{n_0} + o(1).
\end{equation}
To finish the proof of Theorem \ref{thm:one-sided}, let $n$ be in a
middle of the last nonnegative block $\{n_k,\ldots n_{l-1}\}$.  
Then as above, we have
\begin{equation}
\label{eq:logn-diff}
\log { \cal K}^{{\cal P}^+}_{n}  - \log { \cal K}^{{\cal P}^+}_{n_k}
\ge 
\frac{cn}{2}\bar x_n^2 - \frac{c}{n_k} - \frac{c}{2}\log\frac{n}{n_k}.
\end{equation}
Adding 
(\ref{eq:logk-diff}) and
(\ref{eq:logn-diff}) we have
\[
\log { \cal K}^{{\cal P}^+}_{n} - \log { \cal K}^{{\cal P}^+}_{n_0} 
\ge \frac{cn}{2}\bar x_n^2 - \frac{c}{2}\log n - \frac{c}{4}\log n_{k} 
+ \frac{3c}{4}\log n_{0} + o(1).
\]
Noting that $n_k = O(n)$, we derive 
\[
\log { \cal K}^{{\cal P}^+}_{n} 
\ge \frac{cn}{2}\bar x_n^2 - \frac{3c}{4}\log n + O(1) 
= \frac{c}{2}\log n \left(\frac{n\bar x_n^2}{\log n} - \frac{3}{2}\right) 
+ O(1). 
\]
This completes the proof of Theorem \ref{thm:one-sided}.

\section{Multivariate linear protocol}
\label{sec:linear-protocol}

In this section we generalize Theorem \ref{thm:main} to multivariate
linear protocols.  See Section 3 of Vovk, Nouretdinov, Takemura and
Shafer (2005) \cite{wp10} and Section 6 of Takemura and Suzuki (2005)
\cite{takemura-suzuki} for discussions of linear protocols. Since the
following generalization works for any dimension, including the case
of infinite dimension, we assume that Skeptic and Reality choose
elements from a Hilbert space $H$.  The inner product of $x,y\in H$ is
denoted by $x\cdot y$ and the norm of $x \in H$ denoted by $\Vert x
\Vert = (x \cdot x)^{1/2}$.  Actually we do not specifically use
properties of infinite dimensional space and readers may just think of
$H$ as a finite dimensional Euclidean space ${\mathbb R}^m$.  For
example spectral resolution below just corresponds to the spectral
decomposition of a nonnegative definite matrix.

Let ${\cal X}\subset H$ denote the move space of Reality, and assume
that ${\cal X}$ is bounded.  Then by rescaling we can say without loss
of generality that ${\cal X}$ is contained in the unit ball
\[
 {\cal X} \subset \{ x \in H \mid \Vert x \Vert \le 1 \}.
\]
In this case the closed convex hull $\overline{\rm co}({\cal X})$ of
${\cal X}$
is contained in the unit ball. 
As in \cite{takemura-suzuki} we also assume that the origin $0$
belongs to $\overline{\rm co}({\cal X})$.
Note also that the average $\bar x_n$ of
Reality's moves always belongs to $\overline{\rm co}({\cal X})$ and
hence $\Vert \bar x_n\Vert \le 1$.
In order to be clear, we write out our game for
multivariate linear protocol.

\newpage
\medskip
\noindent
\textsc{Bounded Linear Protocol Game in General Dimension}\\
\textbf{Protocol:}

\parshape=6
\IndentI   \WidthI
\IndentI   \WidthI
\IndentII  \WidthII
\IndentII  \WidthII
\IndentII  \WidthII
\IndentI   \WidthI
\noindent
${\cal K}_0 :=1$.\\
FOR  $n=1, 2, \dots$:\\
  Skeptic announces $M_n\in H$.\\
  Reality announces $x_n\in {\cal X}$.\\
  ${\cal K}_n := {\cal K}_{n-1} + M_n \cdot x_n$.\\
END FOR
\medskip

As a natural multivariate generalization of the strategy ${\cal P}^c$ 
given by (\ref{eq:our-strategy}), we consider the strategy 
${\cal P} = {\cal P}^A$ 
\begin{align}
\label{eq:A-strategy}
M_n = A\bar{x}_{n-1}{\cal K}_{n-1},
\end{align}
where $A$ is a self-adjoint operator in $H$.  
Then $A$ has the spectral resolution 
\begin{align}
A = \int_{-\infty}^{\infty} \lambda E(d\lambda),
\end{align}
where $E$ denotes the real spectral measure of $A$, or the 
resolution of the identity corresponding to $A$. 
Let $\sigma(A)$ denote the spectrum of $A$ (i.e.\ the support of $E$)
and let
\[
c_0 =\inf \{\lambda\ |\ \lambda \in \sigma(A)  \},\qquad c_1 = 
\sup \{\lambda\ |\ \lambda \in\sigma(A)  \}.
\]
In the finite dimensional case, $c_0$ and $c_1$  correspond to the
smallest and the largest eigenvalue of the matrix $A$, respectively.

Now we have the following generalization of Theorem
\ref{thm:main}.

\begin{theorem}  \label{thm:hilbert} 
In the bounded linear protocol game in general dimension, 
  if Skeptic uses the strategy (\ref{eq:A-strategy}) with
  $0 < c_0 \le c_1 \le 1/2$.
Then
  $\limsup_n {\cal K}_n =\infty$ for each path $\xi=x_1 x_2 \dots$ of
  Reality's moves such that
\[
\limsup_n \frac{\sqrt{n}\Vert\bar x_n\Vert}{\sqrt{\log n}} > 
\sqrt{\frac{c_1}{c_0}}.
\]
\end{theorem}

\proof
In the expression 
\begin{align}
{\cal K}_n = \prod_{i=1}^n (1 + A\bar{x}_{i-1} \cdot x_i),
\end{align}
we have
\begin{align}
A\bar{x}_{i-1} \cdot x_i = \int_{c_0}^{c_1} 
\lambda (E(d\lambda)\bar{x}_{i-1}
\cdot x_i) = \bar{y}_{i-1} \cdot y_i,
\end{align}
with
\[
\bar{y}_{i-1} = \int_{c_0 }^{c_1}
\sqrt{\lambda}E(d\lambda)\bar{x}_{i-1},\ \ 
y_i = \int_{c_0 }^{c_1}
\sqrt{\lambda}E(d\lambda)x_i.
\]

By the Schwarz's inequality,
\[
|A\bar{x}_{n-1} \cdot x_i| \le 
\Vert \bar{y}_{i-1} \Vert \Vert y_i \Vert 
\le c_1 \Vert \bar{x}_{i-1} \Vert \Vert x_i \Vert \le \frac{1}{2}.
\]
Hence as in (\ref{eq:logkn}), we can bound $\log {\cal K}_n$ from below as 
\begin{align}
\label{eq:mlogkn}
\log {\cal K}_n \ge 
\sum_{i=1}^n \bar y_{i-1}\cdot y_i - c_1 \sum_{i=1}^n \Vert \bar
y_{i-1} \Vert ^2 .
\end{align}
The first term also can be expressed as in (\ref{eq:sum-x-bar2}), 
and it is bounded from below as follows.
\begin{align}
\label{eq:msum-y-bar2}
\sum_{i=1}^n \bar y_{i-1}\cdot y_i 
&= \frac{1}{2} \sum_{i=2}^n \frac{i}{i-1} \Vert \bar  y_i \Vert^2  
+ \frac{n}{2} \Vert \bar y_n \Vert^2
- \frac{1}{2}\Bigl(\Vert y_1 \Vert^2 + 
\sum_{i=2}^n \frac{1}{i-1} \Vert y_i \Vert^2\Bigr) \nonumber \\
&\ge 
\frac{1}{2} \sum_{i=1}^n \Vert \bar  y_i \Vert^2  + 
\frac{n}{2} \Vert \bar y_n\Vert^2
- \frac{c_1}{2}(3 + \log(n-1)) .
\end{align}
Combining (\ref{eq:mlogkn}) and (\ref{eq:msum-y-bar2}), we have
\begin{align}
\label{eq:mlogkn2}
\log {\cal K}_n &\ge
\Bigl(\frac{1}{2} - c_1 \Bigr) \sum_{i=1}^n \Vert \bar y_{i-1}\Vert^2 + 
\frac{n}{2} \Vert \bar y_n\Vert^2 - \frac{c_1}{2}(3 + \log(n-1)) 
\nonumber \\
&\ge \frac{n}{2} \Vert \bar y_n\Vert^2 - \frac{c_1}{2}(3 + \log(n-1))
\nonumber \\
&\ge
\frac{n}{2}\Vert \bar y_n\Vert^2 - \frac{c_1}{2}\log n - \frac{3}{2} c_1
\nonumber \\
&= \frac{c_1}{2} \log n \left( \frac{n \Vert \bar y_n\Vert^2}
{c_1 \log n} -1 \right) - \frac{3}{2}c_1.
\end{align}
It follows that if 
$\limsup_n \sqrt{n}\Vert \bar y_n \Vert/\sqrt{c_1 \log n} > 1$,  
then $\limsup_n
\log {\cal K}_n = +\infty$. 
Now the theorem follows from 
$c_0 \Vert\bar x_n \Vert^2 \le \Vert \bar y_n \Vert^2$.
\qed

\bigskip
Note that 
\[
c_0 \Vert\bar x_n \Vert^2 \le \Vert \bar y_n \Vert^2 \le c_1 \Vert \bar x_n \Vert^2 
\]
and the equalities hold if and only if $A$ is a scalar multiplication 
operator 
$
A = \int_{-\infty}^{\infty} c E(d\lambda).
$

\begin{remark}
Suppose that $\{A_m\}$ is a sequence of positive definite degenerate 
self-adjoint operators with finite dimensional ranges 
$R_{A_m}(H) \varsubsetneq R_{A_{m+1}}(H) \cdots $, and with the supports 
(ranges of eigenvalues) $\sigma(A_m) \varsubsetneq \sigma(A_{m+1}) \cdots 
\subset (0, 1/2]$. Also suppose that $A_\infty$ is a compact operator 
with infinite dimensional range $R_{A_\infty}(H) \subset H$, 
and $A_\infty$ is obtained in the limit $A_m \to A_\infty,\ m \to \infty$. 
Then $c_0(A_m) \to c_0(A_\infty) = 0,\ m \to \infty$, so that in Theorem 
\ref{thm:hilbert}, $\sqrt{c_1(A_m)/c_0(A_m)} \to 
\sqrt{c_1(A_\infty)/c_0(A_\infty)} = \infty$, implying that the 
strategy ${\cal P^{A_\infty}}$ cannot weakly force SLLN with any rate. 
This phenomenon reflects one feature that the dimension of Skeptic's 
move space is related to the effective weakly force of his strategy. 
It is a subject we will treat in the forthcoming paper.
\end{remark}

\appendix
\section{Summary of mixture $\epsilon$-strategy in Chapter 3 of Shafer
  and Vovk (2001)}

Here we summarize the mixture $\epsilon$-strategy in Chapter 3 of Shafer
and Vovk (2001) for the bounded forecasting game in such a way that its
resource requirement (computational time and memory) becomes explicit.

For a single fixed $\epsilon$-strategy ${\cal P}^\epsilon$, which sets $M_n =
\epsilon {\cal K}_{n-1}$,   the capital process is given as
$
{\cal K}_n^{{\cal P}^\epsilon} (x_1 \dots  x_n)=\prod_{i=1}^n (1 +
\epsilon x_i).
$
Shafer and Vovk combine ${\cal P}^\epsilon$ for many values of
$\epsilon$.
Let $ 1/2 \ge \epsilon_1 > \epsilon_2 > \dots > 0$ be a sequence of
positive numbers  converging to $0$.  Write $\epsilon_{-k}=-\epsilon_k$, $k
=1,2,\dots$.  Furthermore let 
$\{ p_i \}_{i=0,\pm 1, \pm 2,\dots}$ 
be a probability distribution on the set of integers
$\mathbb{Z}$, such that $p_0=0$ and $p_{-k}=p_k$ (symmetric).  Actually $p_k$ is the
initial amount (out of 1 dollar) put into the account with the strategy
${\cal P}^{\epsilon_k}$, $k=\pm 1, \pm 2, \dots$.   For example
we could take
\begin{equation}
\label{eq:uniform-like}
\epsilon_k = \frac{1}{2^{|k|}}, \quad
p_k = \frac{1}{2^{|k|+1}}, \ k=\pm 1, \pm 2, \dots, 
\end{equation}
as is done in Chapter 3 of Shafer and Vovk (2001).  However it is more
convenient here to leave $\epsilon_k$ and $p_k$ to be general.
Then the mixture 
strategy weakly forcing  SLLN is given by the weighted average 
of the strategies ${\cal P}^{\epsilon_k}$, $k=\pm 1, \pm 2, \dots$, with
respect to the weights $\{ p_k\}_{k=0,\pm 1, \pm 2,\dots}$, namely
$$
{\cal P}^* = \sum_{k=-\infty}^\infty p_k {\cal P}^{\epsilon_k}.
$$

Consider the value $M_n^*$ of   ${\cal P}^*$.  It is written as
\[
M_n^* = \sum_{k=-\infty}^\infty p_k \epsilon_k \prod_{i=1}^{n-1} (1+\epsilon_k x_i).
\]
Now introduce the elementary symmetric functions $e_{n,0}, e_{n,1}, \dots,
e_{n,n}$ of the numbers $x_1,\dots, x_n$ as
\[
e_{n,0}= 1 , \ e_{n,1}= \sum_{i=1}^n x_i , \ 
e_{n,2} = \sum_{1 \le i < j\le n} x_i x_j, \ \dots, \ 
e_{n,n}= \prod_{i=1}^n x_i.
\]
Recall that the set of values $\{x_1, \dots,x_n \}$ and
the set of values $\{e_{n,1}, \dots,e_{n,n} \}$ are in one-to-one 
relationship.
Therefore computing all the values of the elementary symmetric functions
of $\{x_1, \dots,x_n \}$ is computationally equivalent to keep all the
values of $\{x_1, \dots,x_n \}$.
Using elementary symmetric functions, $M_n^*$ is written as
\begin{equation}
\label{eq:mnstar1}
M_n^* = \sum_{k=-\infty}^\infty p_k \sum_{i=0}^{n-1} \epsilon_k^{i+1} e_{n-1,i}
= \sum_{i=0}^{n-1} e_{n-1,i} \sum_{k=-\infty}^\infty p_k \epsilon_k^{i+1}.
\end{equation}

Let $F$  denote the discrete probability distribution on $[-1/2,1/2]$ with 
$$
F(\{\epsilon_k\})=p_k, \ k=\pm 1, \pm 2, \dots ,
$$
and let 
$$
\mu_m =  E_F(X^m) = \sum_{k=-\infty}^\infty p_k \epsilon_k^m =
\int_{-1/2}^{1/2} x^m F(dx)
$$
denote the $m$-th moment of $F$.  Then (\ref{eq:mnstar1}) is written
as 
\begin{equation}
\label{eq:mnstar2}
M_n^* = \sum_{i=0}^{n-1}  \mu_{i+1} e_{n-1,i} .
\end{equation}

At this point it becomes clear that we can take
an arbitrary probability distribution $F$ on $[-1/2, 1/2]$ provided that
$F$ is not degenerate at 0 and 0 is the point of support, i.e., for
all $\epsilon>0$
\[
F(\epsilon)-F(-\epsilon)>0.
\]
$M_n^*$ in (\ref{eq:mnstar2}) for such an $F$ is a strategy weakly forcing
SLLN.  Now the simplest $F$ seems to be the uniform distribution on
$[-1/2, 1/2]$.  Actually (\ref{eq:uniform-like}) is in a sense close to
the uniform distribution. Then
\[
\mu_m = \int_{-1/2}^{1/2} x^m dx = 
 \begin{cases}
 \frac{1}{m+1} \frac{1}{2^m} & m:\mbox{even} \cr
 0 & m:\mbox{odd}.
\end{cases}  
\]
In this case there is virtually no computational resource
is needed to compute $\mu_m$, since it is explicitly given.  Furthermore
the strategy weakly forcing SLLN based on the uniform distribution is
given by 
\[
M_n^* = \sum_{i=0\atop i:\mbox{\scriptsize odd}}^{n-1}
\frac{1}{i+2} \frac{1}{2^{i+1}} e_{n-1,i} .  
\]
In order to compute the $M_n^*$ we need all the values of the elementary
symmetric functions $e_{n-1,i}$, $i=1,\dots,n-1$, which is equivalent to
keeping  all the values of $\{x_1, \dots,x_{n-1} \}$ as discussed above.
Therefore the mixture $\epsilon$-strategy needs a memory of size proportional
to $n$ at time $n$, whereas our strategy in (\ref{eq:our-strategy}) only
needs to keep track of the values of $\bar x_{n-1}$ and ${\cal K}_{n-1}$.

\end{document}